\documentclass[11pt]{amsart}
\usepackage{graphicx}
\usepackage{amssymb}
\usepackage{epic}

\newtheorem{theorem}{Theorem}
\newtheorem{corollary}[theorem]{Corollary}
\newtheorem{prop}[theorem]{Proposition}
\newtheorem{lemma}[theorem]{Lemma}
\theoremstyle{definition}

\newtheorem{remark}{Remark}

\newcommand{\idiot}[1]{\vspace{5 mm}\par \noindent
\marginpar{\textsc{Note}}
\framebox{\begin{minipage}[c]{0.95 \textwidth}
#1 \end{minipage}}\vspace{5 mm}\par}

\renewcommand{\idiot}[1]{}

\def\NCSI{{NCSym}}
\def\NCQ{{NCQSym}}

\def\QQ{{\mathbb Q}}

\def\NN{{\mathbb N}}

\def\Sym{{Sym}}

\def\la{{\lambda}}

\def\m{{\bf m}}
\def\M{{\bf M}}
\def\w{{\bf w}}
\def\W{{\bf W}}
\def\V{{\bf V}}
\def\q{{\bf q}}
\def\Q{{\bf Q}}
\def\p{{\bf p}}

\def\mleq{{\leq_\ast}}
\def\covers{{\lessdot}}
\def\T{{\mathcal T}}
\def\vbar{\sep \, }
\def\shuf{{\sqcup\!\sqcup}}
\def\sep{{,}}
\def\set#1{\{#1\}}

\title[non-commutative symmetric and quasi-symmetric]{The Hopf algebras 
of  symmetric functions\\
and quasisymmetric functions in non-commutative\\variables are free and cofree}
\author{N. Bergeron}\address[Nantel Bergeron]{Department of Mathematics and Statistics\\ York  University\\       To\-ron\-to, Ontario M3J 1P3\\ CANADA} \email{bergeron@mathstat.yorku.ca}  \urladdr{http://www.math.yorku.ca/bergeron}

\author{M. Zabrocki} \address[Mike Zabrocki]{Department of Mathematics and Statistics\\ York  University\\       To\-ron\-to, Ontario M3J 1P3\\ CANADA} \email{zabrocki@mathstat.yorku.ca}  \urladdr{http://www.math.yorku.ca/\~{}zabrocki}

\date{\today}
 \thanks{This work is supported in part by CRC and NSERC.}

\begin{document}

\begin{abstract}
We uncover the structure of the space
of symmetric functions in non-commutative variables by showing
that the underlined Hopf algebra is both free and co-free.  We also introduce
the Hopf algebra of quasi-symmetric functions in non-commutative
variables and define the product and coproduct on the monomial
basis of this space and show that this Hopf algebra is
free and cofree.  In the process of looking for bases which
generate the space we define  orders on the set partitions
and set compositions which allow us to define bases which
have simple and natural rules for the product of basis elements.
\end{abstract}

\maketitle

\begin{section}{Introduction}

Recent progress has been made on unraveling the structure
of the algebra of symmetric functions in noncommutative variables.
This algebra was originally looked at by Wolf \cite{W} showing that the
algebra is freely generated by elements of the monomial basis.
In the paper by Rosas and Sagan \cite{RS}, bases analogous to
those found in the commutative algebra of symmetric functions
were defined showing that the algebra has similar features to
the algebra of symmetric functions (in commuting variables).  
In \cite{BRRZ}, along with Rosas and Reutenauer,
we introduced a Hopf algebra structure and computed some of the
structure of the coinvariants in the non-commutative polynomials.
In \cite{BHRZ}, along with Rosas and Hohlweg, we found a connection
between the representations of the partition algebra and the dual
of the algebra symmetric functions in non-commutative variables.  
In \cite{HNT} it is remarked that this graded dual can be realized as
a subalgebra of a commutative Hopf algebra of endofunctions and the
techniques there can be used to derive independent proofs that this Hopf
algebra is both free and cofree.
Aguiar and Mahajan \cite{AM} are currently writing a monograph
on combinatorial Hopf algebras and most of these results are 
covered with a geometric perspective.

In this paper we continue to uncover the structure of the space
of symmetric functions in non-commutative variables by showing
that the Hopf algebra is both free and co-free by exhibiting bases
which generate the algebra and co-algebra.  Using similar methods,
we also introduce
the Hopf algebra of quasi-symmetric functions in non-commutative
variables and define the product and coproduct on the monomial
basis of this space and show that this Hopf algebra is
free and cofree.  The quasi-symmetric functions in non-commutative variables
appears in the thesis of 
F. Hivert \cite{H} as an example of combinatorial Hopf algebras
that arise in relation to symmetrization operations.  A personal communication
by J-C. Novelli and J-Y. Thibon informs us \cite{NT} that this Hopf algebra
is bidendriform and can be shown by consequence to be self-dual, free and co-free.

In the process of looking for bases which
generate the space we define  orders on the set partitions
and set compositions which allow us to define bases which
have simple and natural rules for the product of elements.
These orders have interesting properties  in their own
right and are closely related to the structure of the Hopf algebra.

This article is divided into \ref{numsections} sections and an
appendix.  
Sections \ref{sec:notationncsym} and \ref{sec:notationncqsym}
introduce notation and the algebras
of symmetric functions and quasi-symmetric functions in
non-commuting variables and their graded duals (respectively
$NCSym$, $NCQSym$, $NCSym^\ast$ and $NCQSym^\ast$).  
In addition, we
use these sections to quickly run through combinatorial enumerative
results which we use to understand the indexing sets
for the generators of the algebras and their graded duals.

In section \ref{sec:ncsymfree}
we show that $\NCSI$ is freely generated by the analogue of the
power basis defined in \cite{RS} in $\NCSI$.  We also define a second
order on set partitions which defines a multiplicative basis
which does not seem to have a natural anlogue in the algebra of
symmetric functions.  This basis can be used to show that the matrix
of $[ \delta_{A \geq B} ]_{A,B \vdash [n]}$ (where $\leq$ is the refinement
order on set partitions) factors in a very
natural way over the integers.  In section \ref{sec:ncsymcofree}
we show that the graded dual of the algebra $\NCSI$ is freely generated
by a subset of the basis elements which are dual to the monomial basis.
A combinatorial interpretation for the indexing sets follows
from the enumerative results that were found in section \ref{sec:notationncsym}.

In the last two sections we find similar results for
$NCQSym$ and its graded dual.  We introduce two new orders on
set partitions which are used to define natural bases on these
spaces.  These bases have very elegant and simple expressions
for the product of two basis elements and are likely to
play an important role in the structure of these algebras.  The
appendix holds images of one of the  orders on 
set compositions defined in section \ref{sec:ncqsymcofree}.
This order can be generalized to any Coxeter system, 
where one finds that the intervals of the poset are
Cohen-Macauley and several other 
interesting properties. With C. Hohlweg,  this is shown in
a recent paper \cite{BHZ}.

\end{section}
\begin{section}{Set partitions and $NCSym$}\label{sec:notationncsym}

\begin{subsection}{Set Partitions}
Throughout this article, the notation $[n]$  represent  the set 
$\{ 1,2, 3, \ldots, n \}$.  We
say then that $A \subseteq 2^{[n]} \backslash \{\emptyset\}$ with 
$A = \{ A_1, A_2, \ldots, A_k \}$ is a set partition of $n$
(denoted $A \vdash [n]$) if $A_i \cap A_j = \emptyset$
for $i \neq j$ and $A_1 \cup A_2 \cup \cdots \cup A_k = [n]$.  For $A \subseteq 
2^{[n]} \backslash \{\emptyset\}$ with $A_i \cap A_j = \emptyset$ for $i \neq j$
let $st(A)$ represent the set partition formed by lowering the entries of $A$
keeping the values in relative order. The
subsets $A_i$ are referred to as the parts of the set partition.  In this
case we  say the size of the set partition is $|A| = n$ and the length
of $A$ is $\ell(A) = k$ referring to the number of parts of the set partition.  The
number of set partitions of size $n$ are given by the Bell numbers $B_n$ (\cite{OLEIS}
sequence A000110),
and the number of set partitions of length $k$ are the Stirling numbers of
the second kind $S_{n,k}$ (\cite{OLEIS} sequence A008277).  These numbers are given 
by the following recursive formulas:
$$ B_n = \sum_{i=0}^{n-1} { n-1 \choose i } B_{n-1-i} \hbox{ with } B_{0} = B_{1} = 1$$
$$ S_{n,k} = k S_{n-1,k}+S_{n-1,k-1} \hbox{ with }S_{n,1} = S_{n,n} = 1$$

The set partitions of $n$ are ordered by refinement order where
$A \leq B$ if and only if for each $i$, $A_i \subseteq B_j$ for some $j$.
Under this natural order the set $\{ A : A \vdash [n] \}$ is a lattice
where we denote
the g.l.b. of $A$ and $B$ as
$A \wedge B = \{ A_i \cap B_j : 1\leq i \leq \ell(A), 1 \leq j \leq \ell(B) \}$
while $A \vee B$  represent the l.u.b. of $A$ and $B$.

We use the symbol $|$ to represent an
associative operation on set partitions $A\vdash [n]$ and
$B\vdash[k]$.
Let $A|B \vdash [n+k]$ which represents 
$\{ A_1, A_2, \ldots, A_{\ell(A)}, B_1 + n, B_2 + n,
\ldots, B_{\ell(B)}+n \}$ where $B_i + n$ are the entries
of $B_i$ with $n$ added to each. 

A set partition $A \vdash [n]$ is called {\it splittable} if there is a number $r< n$
such that $A_{i_1} \cup A_{i_2} \cup \cdots \cup A_{i_j} = [r]$ for some proper 
subset of the parts of $A$ and $A$ is called {\it atomic} if it is non-empty and
not splittable.  Alternatively we have that $A$ is splitable if
there exists non-empty set partitions $B \vdash [k]$ and $C\vdash [n-k]$ such that
$A = B|C$. For a set partition $A$, we define the split of $A$ to be 
$A^! = ( A^{(1)}, A^{(2)}, \ldots, A^{(d)})$ where each $A^{(i)}$ is an atomic
set partition of size $\alpha_i > 0$ and $A = A^{(1)} | A^{(2)} | \cdots 
| A^{(d)}$. For example, for $A=\set{ \set{1,3},\set{2},\set{4},\set{5,8},\set{6,7}}$ we have $A^!=\big(\set{ \set{1,3},\set{2}},\,\,\set{\set{1}},\,\,\set{\set{1,4},\set{2,3}}\big)$.
By definition we have that $A$ is atomic if and only if $A^! = (A)$.

The bijection between $A$ and $A^!= ( A^{(1)}, A^{(2)}, \ldots, A^{(d)})$ has the property
that $\ell(A) = \ell(A^{(1)})+\ell(A^{(2)}) + \cdots + \ell(A^{(d)})$ and that
implies the following combinatorial result. It is a hint of the type of 
algebraic structure we  find when we consider the symmetric
functions in non-commutative variables.

\begin{prop}
Let $a_{k,i}$ represent the number of atomic set partitions of size $k$
and of length $i$ (\cite{OLEIS} sequence A087903), then we have
\begin{equation}
\frac{1}{1-\sum_{k\geq1} \sum_{i=1}^{k-1} a_{k,i} t^i q^k} 
= \sum_{n\geq0}\sum_{r\geq1} S_{n,r} t^r q^n
\end{equation}
\end{prop}

By setting $t=1$ in the expression above we also have the following corollary.
\begin{corollary}
Let $a_k$ represent the number of atomic set partitions of size $k$
(\cite{OLEIS} sequence A074664), then
\begin{equation}
\frac{1}{1-\sum_{k\geq1} a_{k} q^k} = \sum_{n\geq0} B_n q^n
\end{equation}
\end{corollary}
\end{subsection}
\begin{subsection}{Lyndon words}
Let $x_1, x_2, x_3, \ldots$ be a totally ordered alphabet.  A word in
this alphabet $w = x_{i_1} x_{i_2} \cdots x_{i_{|w|}}$
is called Lyndon if it is lexicographically strictly smaller than all of the cyclic shifts
of $w$.

We have the following proposition about a decomposition of words
into Lyndon words.  This is exercise 7.89.d in \cite{S}:
\begin{prop} (\cite{R} (7.4.1)) \label{Ldecomp}
For every word $w$ in a totally ordered alphabet, there is a
unique decomposition of $w$ as a concatenation of 
Lyndon words $w^{(i)}$ such that
$w = w^{(1)} w^{(2)} \cdots w^{(d)}$ with
$w^{(i)} \geq_{lex} w^{(i+1)}$. 
\end{prop}

Recall that the size of the set partition is denoted by $|A|$ and the
number of parts is $\ell(A)$.
We  consider a total order on set partitions such that $A <_T B$ if
$|A| < |B|$, or $|A|=|B|$ and $\ell(A) < \ell(B)$, or $|A| = |B|$ and $\ell(A)=\ell(B)$
and $w(A) <_{lex} w(B)$ where $w(A)$ is a word with $w(A)_i = j$ if $i \in A^{(j)}$
where $A^!= ( A^{(1)}, A^{(2)}, \ldots, A^{(d)})$. 

The following purely combinatorial result will be used in our discussion
of the primitives of the Hopf algebra of symmetric functions in non-commutative variables.
We say that  $A$, a set partition of $[n]$, is {\it Lyndon} if $A^!$ is Lyndon in lexicographic
order
in the atomic set partitions which are ordered using the total order on
all set partitions listed above.  Remark that every atomic set partition is Lyndon
since $A^!$ has length $1$.

Now by viewing $A^!$ as a word in the alphabet of atomic set partitions 
we may use Proposition \ref{Ldecomp} and the bijection between $A$ and $A^!$
to arrive at the following two corollaries.

\begin{prop}
Let $b_{k,i}$ be the number of Lyndon set partitions of size $k$
and length $i$ (\cite{OLEIS} sequence A112340).  We have
\begin{equation}
\prod_{k \geq 1} \prod_{i \geq 1} \frac{1}{(1-t^i q^k)^{b_{k,i}}}
= \sum_{n\geq0}\sum_{r\geq1} S_{n,r} t^r q^n.
\end{equation}
\end{prop}

\begin{corollary}
Let $b_{k}$ be the number of Lyndon set partitions of size $k$ 
(\cite{OLEIS} sequence A085686).
\begin{equation}
\prod_{k \geq 1} \frac{1}{(1-q^k)^{b_{k}}}
= \sum_{n\geq0}B_n q^n.
\end{equation}
\end{corollary}

\end{subsection}

\begin{subsection}{The Hopf algebra $\NCSI$ and the graded dual}
We refer the reader to \cite{RS} for some of the structure
of the algebra of the symmetric functions in non-commuting variables and
to \cite{BRRZ} for the motivation for the following definition of the Hopf algebra.

Define $\NCSI = \bigoplus_{n \geq 0} {\mathcal L} \{ \m_A : A \vdash [n] \}$
as a graded vector space.  $\NCSI$ is endowed with the following
product and coproduct which makes it a Hopf algebra (see \cite{BRRZ} Theorem 5).
For $A\vdash [n]$ and $B\vdash[k]$, set
\begin{equation}\label{proddef}
\m_A \m_B = \sum_{C \wedge ([n] | [k])=A|B} \m_C.
\end{equation}

The coproduct is defined by
\begin{equation*}
\Delta^{\NCSI}( \m_A ) = \sum_{S \subseteq [\ell(A)]} \m_{st(A_S)} \otimes \m_{st(A_{S^{c}})}
\end{equation*}
where for each subset $S = \{ i_1, i_2, \ldots, i_k \}\subseteq [\ell(A)]$,
we denote $A_S = \{ A_{i_1}, A_{i_2},
\ldots, A_{i_k}\}$.

Using this product and coproduct $\NCSI$ is a Hopf algebra with the unit, counit
defined in the usual manner and the antipode is determined from the graded
bialgebra structure.  This algebra is named the symmetric functions in non-commutative
variables because of the following result:

\begin{prop} (\cite{BRRZ} Corollary 2) 
Define $\NCSI^n \subseteq \QQ \langle x_1, x_2, \ldots, x_n \rangle$
as the linear span of the elements 
\begin{equation*}\m_A[X_n]
=\sum_{\alpha} x_{\alpha_1} x_{\alpha_2} \cdots x_{\alpha_{|A|}}
\hbox{ where }\alpha_r = \alpha_s\hbox{ iff }r,s \in A_i\hbox{ for some }i.
\end{equation*}
The map $\phi_n: \NCSI \rightarrow \NCSI^n$ defined by $\phi_n(\m_A) = \m_A[X_n]$
is an algebra morphism.
\end{prop}

The elements $\m_A[X_n]$ are invariant under the action of the symmetric group
$S_n$ which permutes the variables and may be considered a non-commutative analogue of the
monomial symmetric polynomials, hence we refer to the Hopf algebra $\NCSI$
as the symmetric functions in non-commuting variables (even though we do not
reference the variables when we consider the algebra as the linear span of 
elements $\m_A$).

The graded dual Hopf algebra is denoted $\NCSI^\ast = \bigoplus_{n\geq0}
{\mathcal L} \{ \w_A : A \vdash [n] \}$ where the basis $\w_A$ are the dual elements
to the basis $\{ \m_A \}$ in the sense that there exists a pairing between
$\NCSI$ and $\NCSI^\ast$ with $\left[ \m_A, \w_B \right] = \delta_{AB}$.
The product and coproduct are defined on this space so that 
\begin{equation} \label{dualproddef}
\left[ \m_C, \w_A \w_B \right]  = \left[ \Delta^\NCSI( \m_C ) , \w_A \otimes \w_B \right]
\end{equation}
and 
\begin{equation} \label{dualcoproddef}
\left[ \m_B \otimes \m_C, \Delta^{\NCSI^\ast}(\w_A)\right] = 
\left[ \m_B \m_C, \w_A \right].
\end{equation}
  To this end we explicitly define for $A \vdash [n]$
and $B \vdash [m]$,
\begin{equation} \label{dualprod}
\w_A \w_B = \sum_{S \in { [n+m] \choose n}} 
\w_{A \!\uparrow_S \cup B \!\uparrow_{[n+m]\backslash S}}
\end{equation}
where ${ [n+m] \choose n}$ is the set of $n$ element subsets of $[n+m]$ and
$A \!\uparrow_S$ represents raising the entries in $A$ so that they remain
in the same relative order but so that the union of all of the parts is equal to $S$.

\begin{remark}  Because there is a surjective graded map from
$\NCSI$ to $\Sym$ defined as $\chi( \m_A ) = c_\la m_{\la(A)}$ 
(see \cite{RS} and \cite{BRRZ}), this means that there is an injection of $\Sym^\ast$ to
$\NCSI^\ast$ which is
easy to calculate explictly, namely, $\chi^\ast( h_\la) = c_\la \sum_{\la(A) = \la} \w_A$.
This does not immediately associate the elements $\w_A$ as a `homogeneous' basis
and this is a reason we choose to label them with a $\w$ as upside down `m' rather than an `h.'
\end{remark}

It is easy to give an explicit
expression for the coproduct as well.  For $A \vdash [n]$, 
\begin{align} \label{dualcoprod}
\Delta^{\NCSI^\ast}( \w_A ) &= \sum_{k=0}^n \sum_{([k]|[n-k]) \wedge A = (B|C)} \w_B \otimes \w_C\\
&=\sum_{i=0}^{n} \w_{st(A \!\downarrow_{\{1,\ldots,i\}})} \otimes
\w_{st(A \!\downarrow_{\{i+1, \ldots, n\}})}
\end{align}
where in the sum $B \vdash [k]$ and $C \vdash [n-k]$ and the notation
$A \!\downarrow_S$ is the set partition $A$ restricted to the entries which are in $S$
(throwing away any empty parts).

Since $\NCSI$ is non-commutative and co-commutative, $\NCSI^\ast$ is
commutative and non-co-commutative (as seen in the product
\eqref{dualprod} and coproduct \eqref{dualcoprod} formulas).  We leave as
an exercise to the reader to show that these formulas verify the definitions
of \eqref{dualproddef} and \eqref{dualcoproddef}.

We give the following two examples of the product and co-product to help demonstrate
the formulas above.  The set partitions in the first example below have their parts ordered
so that it is easy to read the set $S \in { [5] \choose 3 }$ from the first 3 numbers that appear in
the set partition.
\begin{align*}
\w_{\{1\vbar23\}} \w_{\{1\vbar2\}} =~&\w_{\{1\vbar23\vbar4\vbar5\}} 
+ \w_{\{1\vbar24\vbar3\vbar5\}} + \w_{\{1\vbar34\vbar2\vbar5\}} 
+ \w_{\{2\vbar34\vbar1\vbar5\}} + \w_{\{1\vbar25\vbar3\vbar4\}}\\
 &+ \w_{\{1\vbar35\vbar2\vbar4\}} + \w_{\{2\vbar35\vbar1\vbar4\}} 
 + \w_{\{1\vbar45\vbar2\vbar3\}} + \w_{\{2\vbar45\vbar1\vbar3\}} + \w_{\{3\vbar45\vbar1\vbar2\}}\\
 =~& \w_{\{1\vbar23\vbar4\vbar5\}} + \w_{\{1\vbar24\vbar3\vbar5\}}
 + \w_{\{1\vbar25\vbar3\vbar4\}} + 2 \w_{\{1\vbar2\vbar34\vbar5\}} 
 + 2 \w_{\{1\vbar2\vbar35\vbar4\}} + 3 \w_{\{1\vbar2\vbar3\vbar45\}}
 \end{align*}
 
 \begin{align*}
 \Delta^{\NCSI^\ast}( \w_{\{ 13\vbar256\vbar4 \}} ) = 
 &\w_{\{ 13\vbar256\vbar4 \}} \otimes 1 + \w_{\{ 13\vbar25\vbar4 \}} \otimes \w_{\{1\}} 
  + \w_{\{ 13\vbar2\vbar4 \}} \otimes \w_{\{12\}}
  + \w_{\{ 13\vbar2 \}} \otimes \w_{\{1\vbar23\}}
 \\&+ \w_{\{ 1\vbar2 \}} \otimes \w_{\{1\vbar34\vbar2\}} 
  + \w_{\{ 1 \}} \otimes \w_{\{2\vbar145\vbar3\}}  + 1 \otimes \w_{\{13\vbar256\vbar4\}}
 \end{align*}

\end{subsection}
\end{section}

\begin{section}{ $\NCSI$ is free }\label{sec:ncsymfree}
One of the earliest references to the algebra of $\NCSI$ is
in a paper by Wolfe \cite{W} who considered it with the goal of 
finding an analogue of the fundamental theorem of symmetric functions.
In that paper it was proved that $\NCSI$ is freely generated by
a subset of basis elements $\m_A$. 

In \cite{BRRZ} we give a combinatorial description of the
basis elements which generate this algebra and give a
modern restatement of one of her set of generators (the other is similar).
To state the result precisely, define for
set partitions $A \vdash [n]$ and $B \vdash [m]$,
$A \ast B = \{ A_1 \cup (B_1+n), A_2 \cup (B_2+n), \ldots, 
A_{\ell} \cup (B_{\ell}+n) \}$ 
(with conventions that $A_r$ and $B_r$ are empty sets for $r$
greater than the respective lengths).  Then say $A$ is non-splitable
if it cannot be written as $A = B \ast C$ for non-empty set
partitions $B$ and $C$.  Note that non-splitable is not the same as atomic and even
the number of non-splitable elements of size $n$ and length $k$
is not equal to the number of atomic set partitions of size $n$
and length $k$.
\begin{prop} (see \cite{W} and \cite{BRRZ})
The algebra $\NCSI$ is freely generated by the
basis elements $\m_A$ where $A$ is non-splitable.
\end{prop}

Here we give a very short combinatorial proof of the result that
$\NCSI$ is free.

Define a basis, $\{ \p_A \}_A$ of $\NCSI$ as
$$\p_A = \sum_{B \geq A} \m_B.$$
This basis was considered in \cite{RS} as the natural analogue of the
power basis of $\Sym$ in the algebra $\NCSI$.  Indeed, in \cite{BHRZ}
we show that this basis is fundamental to $\NCSI$ for representation theoretic
reasons and we prove that it has the following property.  

\begin{prop} (Lemma 4.1.(i) of \cite{BHRZ}) \label{id1}
$$\p_A \p_B = \p_{A|B}$$
\end{prop}

\idiot{
For sake of completeness we include the proof here.
\begin{proof}
Recall that for $A \vdash [n]$ and
$B \vdash [k]$, we have 
$$\m_A \m_B = \sum_{D \wedge ([n]|[k]) = A|B} \m_D.$$

Now consider
\begin{align}
\p_A \p_B &= \sum_{C \geq A} \sum_{D \geq B} \m_C \m_D \nonumber\\
&= \sum_{C \geq A} \sum_{D \geq B} \sum_{E \wedge ([n]|[k]) = C|D} \m_E \label{eq1}
\end{align}
Notice that we have that $E \geq A|B$ and for each $E \geq A|B$ there is
exactly one term in the triple sum of \eqref{eq1}.  In fact, this implies that the sum
is equal to
$$= \sum_{E \geq A|B} \m_E = \p_{A|B}.$$
\end{proof}}

We can conclude from this result the following corollary.
\begin{corollary}
$\NCSI$ is freely generated by the elements $\p_A$ where $A$ is atomic.
\end{corollary}
\begin{proof}
Since the $\p_A$ are multiplicative, we have that for $A^! = (A^{(1)}, A^{(2)},
\ldots, A^{(k)})$,
$$\p_A = \p_{A^{(1)}} \p_{A^{(2)}} \cdots \p_{A^{(k)}}.$$
Because the $\p_A$ are linearly independent, the $\{ \p_A: A\hbox{ atomic}\}$
must be algebraically independent.
\end{proof}

Before departing from the question of freeness of $\NCSI$ we would
like to introduce another basis which is multiplicative and arises naturally
in the study of this space.  To this end we define the following
 order on the collection of set partitions of size $n$.
Say that $B$ covers $A$ \idiot{(denoted $B \covers A$)} in our  order if and only if 
there are two parts of $A$, $A_i$ and $A_{i'}$ with all elements
in $A_i$ less than all elements in $A_{i'}$ and 
$B = (A \slash \{ A_i, A_{i'} \}) \cup \{ A_i \cup A_{i'}\}$.
Define $A \mleq B$ as the closure of this covering relation.

For $n=3$ and $n=4$ we have the following poset diagrams.
\begin{figure}[htbp] 
   \centering
   \includegraphics[width=16cm]{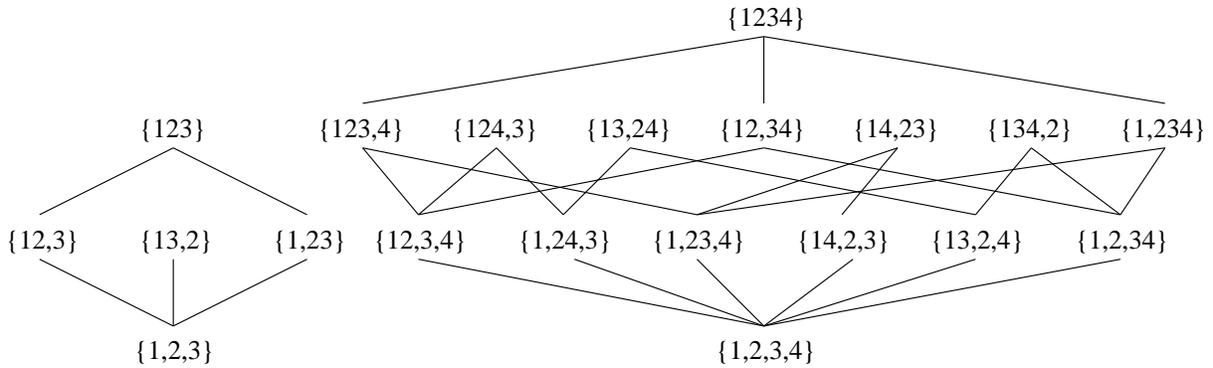} 
   \caption{Hasse diagram of set partitions of size $3$ and $4$ with $\leq_\ast$ order.}
   \label{fig:example}
\end{figure}

Notice that the maximal elements in these two posets are all atomic, although
not all atomic set partitions are maximal (e.g. $\{14\vbar2\vbar3\}$ is atomic 
and not maximal).

\begin{prop} The 
 order $\mleq$ is ranked by $n - \ell(A)$ 
and the poset is Eulerian.  Moreover, fixing a set partition
$A$, the sub-poset consisting
of the elements $\{ B: B \mleq A\}$ with the $\mleq$ order
is a boolean lattice.
\end{prop}

\begin{proof}
In the covering relation defining the $\mleq$ order we have that $B$
is of length $1$ less than the length of $A$ and hence the poset 
is ranked.  The minimal element of this poset is $\{1\vbar2\vbar\cdots\vbar n\}$
and hence the rank of any set partition $A$ in this poset is $n - \ell(A)$. 

Now let ${\mathcal P}_A = \{ B: B \mleq A\}$ be the sub-poset with
maximal element $A$.  For any element $B \in {\mathcal P}_A$, let $\phi(B)$
be a subset of $[n]$ which consists of the elements which are not
maximal in their parts in $B$.  Now if $B$ covers $C \in {\mathcal P}_A$, then
$\phi(C) \cup \{ max(C_i) \} = \phi(B) $ since $B$ is formed 
from $C$ by joining two parts $C_i$ and $C_{i'}$ where all elements
of $C_i$ are less than all elements of $C_{i'}$.  This is the covering
relation for the boolean lattice and hence ${\mathcal P}_A$ is isomorphic to the
boolean lattice on the set $\phi(A)$.

This also implies that any interval of 
the poset is Eulerian since the M\"obius function
at $A$ depends only on the interval ${\mathcal P}_A$. Since this is a boolean
lattice the M\"obius function is equal to $(-1)^{rank(A)}$.
\end{proof}

We observe in this poset that the number of elements covered by $A$
is equal to $n - \ell(A)$.  This is because the number of places that a
part of $A$ can be split is equal to the number of elements of $\phi(A)$.

Define now a new basis of $\NCSI$ by
\begin{equation*}
\q_A = \sum_{A \mleq B} \m_B
\end{equation*}

The property that makes this basis intrinsic to this algebra is
that the elements $\q_A$ where $A$ is atomic are multiplicative, hence
are algebraic
generators of this basis.
We leave the proof of the following theorem to the diligent reader.

\begin{theorem} \label{th:pfree}
For $A \vdash [n]$ where $A^! = (A^{(1)}, A^{(2)}, \ldots, A^{(k)})$ we have
that
\begin{equation*}
\q_A = \q_{A^{(1)}} \q_{A^{(2)}} \cdots \q_{A^{(k)}}.
\end{equation*}
Moreover, the algebra $\NCSI$ is freely generated by the elements
$$\{ \q_A \,|\, A \vdash [n] , n\geq 0, A \hbox{ is atomic}\}.$$
\end{theorem}

\idiot{
\noindent
{\bf Example: }
The atomic set partitions of size less than $5$ are $\{ 1 \}$, $\{12\}$, $\{123\}$, $\{13\vbar2\}$,
$\{1234\}$, $\{124\vbar3\}$, $\{134\vbar2\}$, $\{14\vbar23\}$, $\{13\vbar24\}$, $\{14\vbar2\vbar3\}$.  
For each of these set partitions
except for $B = \{14\vbar2\vbar3\}$ we have that $\q_B = \m_B$ and $\q_{\{14\vbar2\vbar3\}} =
\m_{\{14\vbar2\vbar3\}} + \m_{\{14\vbar23\}}$.  As an example we choose the
set partition $B = \{12\vbar36\vbar4\vbar5\}$ and so $B^! = (\{12\}, \{14\vbar2\vbar3\})$.
Now we calculate
\begin{align*}
\q_{\{12\vbar36\vbar4\vbar5\}} &= \m_{\{12\vbar36\vbar4\vbar5\}} 
+ \m_{\{1236\vbar4\vbar5\}} + \m_{\{124\vbar36\vbar5\}} +
\m_{\{125\vbar36\vbar4\}} + \m_{\{12\vbar36\vbar45\}} + \m_{\{1236\vbar45\}} \\
&= \m_{\{12\}} (\m_{\{14\vbar2\vbar3\}} + \m_{\{14\vbar23\}}) = \q_{\{12\}} \q_{\{14\vbar2\vbar3\}}
\end{align*}
}
\idiot{
Before we proceed with the proof of this theorem we  make some
observations about the relationship between the order and the product.

\begin{lemma} \label{lem:restr}
If $A \mleq B$ then for all subsets $S \subseteq [n]$, 
$st( A\!\downarrow_S ) \mleq st( B\!\downarrow_S)$.
\end{lemma}

\begin{proof}
It suffices to show this for $A \covers B$.  We may thus
assume that $B = A \slash \{ A_i, A_{i'} \} \cup \{ A_i \cup A_{i'} \}$.
If we have that $S \cap A_i = \emptyset$ or $S \cap A_{i'} = \emptyset$,
then $B \!\downarrow_S = A \!\downarrow_S$ which implies that
$st( A \!\downarrow_S ) = st( B \!\downarrow_S)$.

Otherwise we have that $B \!\downarrow_S = A \!\downarrow_S \slash \{
A_i \!\downarrow_S, A_{i'} \!\downarrow_S \} \cup \{ A_i \cup A_{i'} \!\downarrow_S \}$.
Since $st(B\!\downarrow_S) = st( A \!\downarrow_S ) \slash \{ {\bar A}_i , {\bar A}_{i'} \} \cup
\{ {\bar A}_i \cup {\bar A}_{i'} \}$ for sets ${\bar A}_i , {\bar A}_{i'}$ which are parts
of $st( A \!\downarrow_S )$, we have by definition that
$st( A \!\downarrow_S ) \covers st(B\!\downarrow_S)$.
\end{proof}
}

\idiot{
\begin{lemma} \label{lem:splitorder}
Let $A, B \vdash [n]$ and $C, D \vdash [m]$ such that $A \mleq B$ and $C \mleq D$
then $A |C \mleq B |D$.
\end{lemma}
\begin{proof}
Assume there are chains $A \covers B' \covers B'' \covers \cdots \covers B$
and $C \covers D' \covers D'' \covers \cdots \covers D$.  Then
\begin{equation*}
A | C \covers B' |C \covers B'' |C\covers \cdots \covers 
B|C \covers B|D\covers B|D''
\covers \cdots \covers B |D.
\end{equation*}
\end{proof}
}

\idiot{
\begin{proof} (of Theorem \ref{th:pfree})
We proceed by induction on the length of $A^!$.  Assume that for all
set partitions where $A^!$ has length $k-1$ or less that
statement of the theorem holds.  This is vacuously true for the base
case of $k=1$.

Now when $A^! = (A^{(1)}, A^{(2)}, \ldots, A^{(k)})$, take ${\tilde A}$ to
be the set partition with ${\tilde A}^! = (A^{(1)}, A^{(2)}, \ldots, A^{(k-1)})$
and $r = |{\tilde A}| = n - |A^{(k)}|$.
That is, ${\tilde A} = A \!\downarrow_{\{ 1, \ldots, r\} }$ and
$A^{(k)} = st( A \!\downarrow_{\{r+1, \ldots n\}})$.
We have then that
\begin{equation}
\q_{A^{(1)}} \q_{A^{(2)}} \cdots \q_{A^{(k)}} = \q_{\tilde A} \q_{A^{(k)}}\nonumber
= \left(\sum_{{\tilde A} \mleq B} \m_B\right) \left(
\sum_{A^{(k)} \mleq C} \m_C \right).\label{eq:pprod}
\end{equation}

Now consider any set partition $D$ with $A \mleq D$.
We have by Lemma \ref{lem:restr}
that $st( D \!\downarrow_{\{r+1, \ldots, n\}}) = C'$ where $C'$ is some
set partition with $A^{(k)} \mleq C'$
and $st( D \!\downarrow_{\{1, \ldots, r \}}) = B'$ where ${\tilde A} \mleq B'$.
Using \eqref{dualcoproddef} and \eqref{dualcoprod}, the coefficient of 
$\m_D$ in $\m_{B'} \m_{C'}$ is
$[ \m_{B'} \m_{C'}, \w_D] = [\m_{B'} \otimes \m_{C'}, \Delta( \w_D )] = 1$.
Now, $0 < [ \m_{B'} \m_{C'}, \w_D] = [\m_{B'} \otimes \m_{C'}, \Delta( \w_D )]$
implies by equation \eqref{dualcoprod} 
that $D \!\downarrow_{\{1,\ldots,|A|\}} = B'$ and $D \!\downarrow_{\{|A|+1, \ldots,
|D|\}} = C'$.  But then $D$ is the join of $B'$ and $C'$ and hence $D = B'| C'$.
Therefore the only product of $\m_B \m_C$ that $\m_D$ appears in is $B=B'$
and $C=C'$, hence $[ \q_A, \w_D ] = 1$.

Say now that we know
that $[ \q_{\tilde A} \q_{A^{(k)}}, \w_D] >0$ implies $[ \m_B \m_C, \w_D] > 0$
for some $B \vdash [r]$ and $C \vdash [n-r]$ with ${\tilde A} \mleq B$ and
$A^{(k)} \mleq C$.  This implies by Lemma \ref{lem:splitorder} that
 $A =  {\tilde A}  | A^{(k)} 
\mleq B | C \mleq D$.  We have shown that
$[ \q_{\tilde A} \q_{A^{(k)}}, \w_D] = 1$ if and only if $A \mleq D$, therefore
we have $\q_{\tilde A} \q_{A^{(k)}} = \q_A$, completing the proof by induction.

This also shows that the elements $\q_A$, for $A$ atomic, freely generate the
algebra since every basis element may be written as a product of these
elements.
\end{proof}
}

\idiot{
Because we are working with a basis defined by an Eulerian poset
we also have an easy way of defining the dual basis within
$\NCSI^*$.  Define
\begin{equation*}
\q^\ast_A = \sum_{B \mleq A} (-1)^{\ell(B) - \ell(A)} \w_B
\end{equation*}
Because these bases are defined as a sum over elements in
an Eulerian poset we have by M\"obius inversion that
\begin{align*}
\left[ \q_A, \q^\ast_B \right] &= \left[ \sum_{A \mleq A'} \m_{A'},
 \sum_{B' \mleq B} (-1)^{\ell(B') - \ell(B)} \w_{B'} \right]\\
&= \sum_{A \mleq C \mleq B} (-1)^{\ell(C) - \ell(B)} = \delta_{AB}
\end{align*}

The basis $\q_A$ and the order which defines it is also very
interesting because of what it says about the
the refinement order on set partitions.  We observe that

\begin{equation}\label{eq:pinqeq}
\p_A = \q_A + \sum_{B \geq A, A \nleq_\ast B} \q_B.
\end{equation}

For example, observe that
\begin{align*}\p_{\{13\vbar2\vbar4\}} = \m_{\{13\vbar2\vbar4\}}+\m_{\{13\vbar24\}}+
\m_{\{134\vbar2\}}+\m_{\{123\vbar4\}}+\m_{\{1234\}} 
= \q_{\{13\vbar2\vbar4\}}+\q_{\{123\vbar4\}}
\end{align*}
since $\q_{\{13\vbar2\vbar4\}} = \m_{\{13\vbar2\vbar4\}} + \m_{\{134\vbar2\}}
+ \m_{\{13\vbar24\}}$ and $\q_{\{123\vbar4\}} = \m_{\{123\vbar4\}} + \m_{\{1234\}}$.

What equation \eqref{eq:pinqeq} says is that the change of basis matrix between
the $\p_A$ and the $\m_A$ basis (which is the $\zeta$-function for the poset
of the set partitions) factors in a very natural way over the integers into a
product of two unitriangular matrices.
}

\idiot{
\begin{remark} 
The change of basis matrix from $\q$ to $\m$ has entries in $\{ 0,1\}$ and from $\m$ to $\q$
has entries in $\{ -1, 0,1\}$ (the poset is Eulerian).  The change of basis from $\p$ to $\m$
has entries in $\{ 0, 1\}$ and from $\m$ to $\p$ is $(-1)^{\ell(B)-\ell(A)} (\ell(B)-\ell(A))!$.
The change of basis from $\p$ to $\q$ has entries in $\{ 0,1\}$ and experimentally we
observe that from $\q$ to $\p$ has mostly entries in $\{ -1,0,1\}$ (at $n=5$ I found some
values of $2$).  We should be able to
explicitly compute this matrix.
\end{remark}
}

\end{section}
\begin{section}{$\NCSI$ is cofree}\label{sec:ncsymcofree}

We  show in this section that the algebra $\NCSI^\ast$ is generated algebraically
by the elements $\w_A$ such that $A$ is Lyndon (recall that $A$ is Lyndon if $A^!$ is
a Lyndon word in the alphabet of set partitions).  To begin, we consider the shuffle
of two words $u \shuf v$ to be the multi-set of all shuffles of the words $u$ and $v$.
In this definition we  set for $S  \in { [|u|+|v|] \choose |u| }$,
$u \shuf_S v$ to be the word with the $i^{th}$ letter in $u$ if $i \in S$ and
a letter of $v$ if $i \notin S$ (the letters in relative order of the word $u$ and $v$).
We are then considering $u \shuf v$ to be the multi-set of words
$\{ u \shuf_S v : S \in { [|u|+|v|] \choose |u| } \}$.

Given a totally ordered alphabet $X=x_1, x_2,x_3,\ldots$, we denote by $X^*_\shuf$ the shuffle algebra. This is the commutative algebra of all words in the alphabet $X$ with the shuffle product. We will use the following proposition about shuffle algebras.

\begin{prop}\label{ReutShuff} (\cite{R} Theorem 6.1) For any alphabet $X$, the algebra $X^*_\shuf$ is freely generated by the elements of the set of words $\set{u|\,u\,  \hbox{\it  is Lyndon}}$
\end{prop}

As mentioned in the first section, a set partition $A$ corresponds to
a list of atomic set partitions $A^!$ which we may consider as a word
in the alphabet of atomic set partitions.

\begin{lemma}  Let $A \vdash [n]$ and $B \vdash [m]$.
For $S \in {[n+m] \choose n}$ then
$$\ell( (A\!\uparrow_S \cup B \!\uparrow_{S^c})^!) \leq \ell(A^!) + \ell(B^!)$$
with equality if and only if $(A\!\uparrow_S \cup B \!\uparrow_{S^c})^! \in A^! \shuf B^!$.
\end{lemma}

\begin{proof}
Let $A^!=(A^{(1)},A^{(2)},\ldots,A^{(d)})$ and $B^!=(B^{(1)},B^{(2)},\ldots,B^{(p)})$.  For any $S\in  {[n+m] \choose n}$, the atomic parts $A^{(i)}\!\uparrow_S$ [resp. $B^{(j)}\!\uparrow_{S^c}$] of $(A\!\uparrow_S)^!$ [resp. $(B\!\uparrow_{S^c})^!$] can never be split in smaller parts in $(A\!\uparrow_S \cup B \!\uparrow_{S^c})^! $. This implies that the atomic parts of $(A\!\uparrow_S \cup B \!\uparrow_{S^c})^!$ combines parts of $(A\!\uparrow_S)^!$  and $(B\!\uparrow_{S^c})^!$.  These correspond to the parts of $A^!$ and $B^!$. The inequality follows from this fact. For the equality, we remark that this happen if and only if no atomic parts of $A^!$ and $B^!$ are combined, hence $(A\!\uparrow_S \cup B \!\uparrow_{S^c})^!$ is a shuffle of $A^!$ and $B^!$. For the converse we remark that any word in the shuffle $A^! \shuf B^!$ is obtained by a unique choice of $S$.
\end{proof}

Now because $\w_A \w_B = \sum_{S \in {[n+m] \choose n}} \w_{A\!\uparrow_S \cup B \!\uparrow_{S^c}}$
this tells us that the sum may be broken up into terms $\w_C$ where $C^!$ has
maximal length and all other terms.

\begin{corollary}\label{shuffTriang}
\begin{equation*}
\w_A \w_B = \sum_{C : C^! \in A^! \shuf B^!} \w_C + \sum_{D : \ell(D^!) < \ell(A^!)+\ell(B^!)} \w_D
\end{equation*}
\end{corollary}

In order to find algebraic generators, we  find a total order such that $\NCSI^\ast$  is isomorphic to a shuffle algebra.  To this end we remark that given a total order on the atomic set partitions 
we get a lexicographic order $\leq_{lex}$ on words in the atomic set partitions.
Now for $A^{!} = (A^{(1)}, A^{(2)}, \ldots, A^{(k)})$ we will denote the number of atomic
set partitions in $A^!$ by $\ell(A^!)= k$.
We then define
$A <_\T B$ if $\ell(A^!) < \ell(B^!)$, or $\ell(A^!) =\ell(B^!)$ and $A^! <_{lex} B^!.$  Using this order  in Corollary~\ref{shuffTriang} we have that all terms in the second sum 
are strictly smaller then all terms in the first sum. By triangularity, this shows that $\NCSI^\ast$ is isomorphic to the shuffle algebra in the alphabet given by atomic set partitions.
Combining this with Proposition~\ref{ReutShuff}, we have shown the following theorem.

\begin{theorem}  
The algebra $\NCSI^\ast$ is freely generated by the elements
$$\{ \w_B \,|\, B \vdash [n] , n\geq 0, B \hbox{ is Lyndon}\}.$$
\end{theorem}

\end{section}

\begin{section}{Set compositions and $NCQSym$}\label{sec:notationncqsym}
\begin{subsection}{Set compositions }

There is an algebra related to $\NCSI$ which is indexed by set compositions
(ordered set partitions).  A set composition of $n$ is a sequence
$\Phi = ( \Phi_1, \Phi_2, \ldots, \Phi_{\ell(\Phi)})$ where each
$\Phi_i$ is a non empty subset of $[n]$ and $\Phi_i \cap \Phi_j = \emptyset$
if $i \neq j$ (note that we are using $\ell(\Phi)$ to represent the number
of parts or the length of the set composition $\Phi$).
As with set partitions, 
the size of the set composition will be denoted $|\Phi| = |  \Phi_1|+| \Phi_2|
+\cdots+ |\Phi_{\ell(\Phi)}|$.
The number of set compositions of size $n$
and of length $k$ is equal to $k! S_{n,k}$ (\cite{OLEIS} sequence A019538)
the total number of set compositions of $n$ is of course 
$\sum_{k=1}^n k! S_{n,k}$ (\cite{OLEIS} sequence A000670)
are sometimes known as the ordered Bell numbers.

As we used symbols
$A, B, C, D$ to represent set partitions, we will generally use
symbols $\Phi, \Psi, \Pi$ and $\Gamma$ to represent 
set compositions and use the notation $\Phi \models [n]$ to
indicate that $\Phi$ is a set composition of $n$.
In addition, when writing explicit set compositions we will leave
off any $\{ \}$ around the parts and commas removed from 
between the elements of each of the sets for brevity.  
For example, the set composition $(\{3,4\}\vbar\{2\}\vbar\{1,5\})$ will be represented by
$(34\vbar 2\vbar 15)$.

The set compositions of $n$ are endowed with a 
natural  order similar to the refinement order on set 
partitions.  We will say that for $\Phi, \Psi \models [n]$,
$\Phi \leq \Psi$ if for each $i$, $\Phi_i \subseteq \Psi_j$ for
some $j$ and if $\Phi_i \subseteq \Psi_j$, then either
$\Phi_{i+1} \subseteq \Psi_j$ or $\Phi_{i+1} \subseteq \Psi_{j+1}$.

This is a  order with covering relations
$\Phi \covers (\Phi_1, \ldots, \Phi_i \cup \Phi_{i+1}, \Phi_{i+2}, \ldots, \Phi_{\ell(\Phi)})$
for each $1 \leq i < \ell(\Phi)$.  This order does not form a lattice
on the set of set compositions of $n$, but it is a ranked poset with
no single minimal element (the permutations correspond to all
of the minimal elements of this poset) and a maximal element $([n])$.
As an example, we show the diagram of the poset of set compositions
of size 3 defined by this relation here.

\begin{figure}[htbp] 
   \centering
   \includegraphics[width=10.7cm]{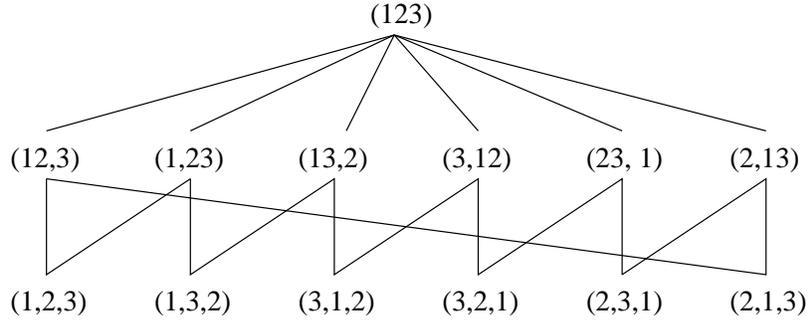} 
   \caption{Hasse diagram of refinement order for set compositions of size 3.}
   \label{fig:example}
\end{figure}
\setlength{\unitlength}{0.00083333in}

Although there is not a lattice structure, there are two operations
similar to the meet and join operations on set partitions
for the poset of set compositions which we shall exploit for notational
purposes.  For two set partitions $\Phi, \Psi \models [n]$, 
there is a well defined element $\Phi \vee \Psi$ which is
the smallest element which is larger than both $\Phi$ and
$\Psi$.  There is not necessarily a single element which
is a greatest element less than both $\Phi$ and $\Psi$, so
we will choose one canonical one as a representative.
Define
\begin{equation}
\Phi \wedge \Psi = (\Phi_1 \cap \Psi_1, \Phi_1 \cap \Psi_2, \ldots, 
\Phi_2 \cap \Psi_1, \Phi_2 \cap \Psi _2
\ldots, \Phi_{\ell(\Phi)} \cap \Psi_{\ell(\Psi)}).
\end{equation}
Note that in general $\Phi \wedge \Psi \neq \Psi \wedge \Phi$
as one would normally have in a lattice structure.  
It will hold that $\Phi \wedge \Psi \leq \Phi$, but in general we
do not have $\Phi \wedge \Psi \leq \Psi$.

The operations which exist on set partitions can be extended
in a natural way to set compositions.  For instance, define
for $\Phi \models[n]$ and $\Psi \models [k]$,
$$\Phi | \Psi = ( \Phi_1, \ldots, \Phi_{\ell(\Phi)}, \Psi_1+n, \ldots,
\Psi_{\ell(\Psi)} + n).$$
It is easily checked that each of the operations $\vee$, $\wedge$ and
$|$ are associative.  Of the three, only the $\vee$ operation is commutative.
Also, for $S\in { [n+m] \choose n}$ we denote by
$\Phi\!\uparrow_S$ the set composition obtained by raising the entries 
in $\Phi$ so that they remain
in the same relative order but the entries are the elements of $S$. 
As well if $\Psi$ is a set composition of an arbitrary finite set of integers,  
$st(\Psi)$ represent the set composition formed by lowering the entries 
of $\Psi$ keeping the values in the same relative order.

With the definition of the $|$ operation we may also define
the concept of `splitable' with respect to this operation.  We
will say that $\Phi \models [n]$ is splitable if there exists non-empy set
compositions $\Psi \models [k]$ and $\Gamma \models [n-k]$
such that $\Phi = \Psi|\Gamma$.  If $\Phi$ is not splitable then
it will be called atomic.   Just as we did for set partitions,
we will define a bijection between set compositions $\Phi$ and
sequences of set compositions
$\Phi^! = (\Phi^{(1)}, \Phi^{(2)}, \ldots, \Phi^{(k)})$
where $\Phi^{(r)}$ are each non-empty atomic set compositions
such that $\Phi = \Phi^{(1)}|\Phi^{(2)}| \cdots |\Phi^{(k)}$.

Since we have that $\ell(\Phi) = \ell(\Phi^{(1)})+\ell(\Phi^{(2)})+
\cdots+\ell(\Phi^{(k)})$ this bijection implies the following enumerative result
regarding set compositions.

\begin{prop}
Let $c_{k,i}$ represent the number of atomic set compositions of size $k$
and of length $i$ (\cite{OLEIS} sequence A109062), then we have
\begin{equation}
\frac{1}{1-\sum_{k\geq1} \sum_{i=1}^{k-1} c_{k,i} t^i q^k} 
= \sum_{n\geq0}\sum_{r\geq1}  r! S_{n,r} t^r q^n
\end{equation}
\end{prop}

There is a natural map from the set compositions to the set partitions
where one `forgets' the order on the parts of the set composition.
We will use the notation $A( \Phi)$ to represent the set partition
$\{ \Phi_1, \Phi_2, \ldots, \Phi_{\ell(\Phi)} \}$.
It follows that this map is compatible with the lattice of set partitions in
the sense that
$A(\Phi|\Psi) = A(\Phi) | A(\Psi)$, $\ell(A(\Phi)) = \ell(\Phi)$, $A(\Phi \vee \Psi) = A(\Phi) \vee
A(\Psi)$, $A(\Phi \wedge \Psi) = A(\Phi) \wedge A(\Psi)$ and
if $\Phi \leq \Psi$ then $A(\Phi) \leq A(\Psi)$ .





In the next section we will define a bialgebra related to $\NCSI$
which is both non-commutative and non-cocommutative (that is, the graded
dual of this algebra will also be non-commutative). 
\end{subsection}

\begin{subsection}{ The Hopf algebra of $\NCQ$ and its graded dual }
As we defined an analogue of the symmetric functions
which are non-commutative, one may also define an analogue of the
quasi-symmetric functions in non-commuting variables.  Recall that
a polynomial is quasi-symmetric if for every increasing sequence
of indices $i_1 < i_2 < \cdots < i_{k}$ and every composition $\alpha
= (\alpha_1, \alpha_2, \cdots, \alpha_k)$ of length $k$ (
with $\alpha_i>0$) we have that the coefficient of the monomial
$x_{i_1}^{\alpha_1} x_{i_2}^{\alpha_2} \cdots x_{i_k}^{\alpha_k}$
is equal to the coefficient of the monomial $x_1^{\alpha_1} x_2^{\alpha_2}
\cdots x_k^{\alpha_k}$.

Now for a sequence $\gamma = (\gamma_{1}, \gamma_2, \ldots, \gamma_n) \in \NN^n$, we let
$\Psi = \Delta( \gamma )$ be a set composition $\Psi \models [n]$ such that
$i \in \Psi_{\#\{ \gamma_r : \gamma_r < \gamma_i \} +1}$.  For example,
$\Delta( 2, 1, 1, 7, 9, 1, 2, 7) = ( 236\vbar17\vbar48\vbar5 )$.

We will define a polynomial in $n$ non-commuting variables
$\{ x_1, x_2, \ldots x_n \}$ to be quasi-symmetric
if for every pair of sequences $\gamma, \tau \in {[n]}^k$
such that $\Delta(\gamma) = \Delta(\tau)$, the coefficient of
$x_{\gamma_1} x_{\gamma_2} \cdots x_{\gamma_k}$ is equal to the 
coefficient of  $x_{\tau_1} x_{\tau_2} \cdots x_{\tau_k}$.


There is a natural basis of the space of quasi-symmetric polynomials
which are similar to monomial symmetric functions and it is this basis
for which we define an analogue.  For an ordered set of variables $X$,
and for a set composition of $n$ which is of length $k$
define
\begin{equation}
\M_{\Phi}[X] = \sum_{\gamma} x_{\gamma_1} x_{\gamma_2} \cdots x_{\gamma_n}
\end{equation}
where the sum is over all sequences $\gamma = (\gamma_1, \gamma_2, \ldots, \gamma_n) \in \NN^n$
such that $\Delta(\gamma)= \Phi$.

Just as we did for the symmetric functions in non-commutative variables we can
define a product and coproduct on this space which endows it with a
Hopf algebra structure.  Without giving reference to the variables
we define an abstract vectors space by
\begin{equation}
\NCQ_n = {\mathcal L} \{ \M_\Phi : \Phi \models [n] \}
\end{equation}
and then set $\NCQ = \bigoplus_{n \geq 0} \NCQ_n$.

The product is defined so that it inherits the
product $\M_\Phi[X] \M_\Psi[X]$ from the space of non-commutative
polynomials.  For $\Phi \models [n]$ and $\Psi \models [k]$
\begin{equation}
\M_\Phi \M_\Psi = \sum_{(([n])|([k])) \wedge \Gamma = \Phi | \Psi} \M_\Gamma.
\end{equation}

In addition we can define a coproduct structure which is
the natural analogue of the coproduct structure on the other
spaces since it follows by essentially replacing one set of
variables by two.  Define the map $\Delta : \NCQ_n \rightarrow
\bigoplus_{k=0}^n \NCQ_k \otimes \NCQ_{n-k}$ by
\begin{equation}
\Delta( \M_\Phi) = \sum_{i=0}^{\ell(\Phi)} \M_{st(\Phi_1,\ldots, \Phi_i)}
\otimes \M_{st(\Phi_{i+1}, \ldots, \Phi_{\ell(\Phi)})}.
\end{equation}
There is a natural embedding of $\NCSI$ to $\NCQ$
just as there is a relationship between the symmetric
functions and the quasi-symmetric functions.  The monomial
basis of symmetric functions of $\NCSI$ are related to
the monomial basis of $\NCQ$ by the map $\theta : \NCSI \rightarrow \NCQ$
which is defined by
\begin{equation}
\theta(\m_A) = \sum_{A = A(\Phi)} \M_\Phi.
\end{equation}

We will also be interested in the graded dual of the
algebra $\NCQ$.  We will define this algebra as the
linear span of the elements $\W_\Phi$ which are
the basis which is dual to the monomial basis $\M_\Phi$.
Set $\NCQ^\ast_n = {\mathcal L}\{ \W_\Phi : \Phi \models [n] \}$
and then $\NCQ^\ast = \bigoplus_{n \geq 0} \NCQ^\ast_n$.  The
pairing between $\NCQ$ and $\NCQ^\ast$ defined by
$[ \M_\Phi, \W_\Psi ] = \delta_{\Phi\Psi}$ defines the product
and coproduct on $\NCQ^\ast$ through the duality relations
\begin{equation}
[ \Delta( \M_\Phi ), \W_\Psi \otimes \W_\Gamma ] = [ \M_\Phi, \W_\Psi \W_\Gamma ]
\end{equation}
and
\begin{equation}
[ \M_\Phi \M_\Psi, \W_\Gamma ] = [ \M_\Phi \otimes \M_\Psi, \Delta^\ast( \W_\Gamma ) ].
\end{equation}

From these two relations it is easy to determine that for set compositions $\Phi \models
[n]$ and $\Psi \models [k]$ we have
\begin{equation} \label{eq:Wprod}
\W_\Phi \W_\Psi = \sum_{S \in {[n+m] \choose n}} \W_{\Phi\!\uparrow_S \cdot 
\Psi\!\uparrow_{S^c}}
\end{equation}
where $\cdot$ represents concatenation. Finally, the coproduct  $\Delta^\ast$ is given explicitly on the dual basis as
\begin{equation}
\Delta^\ast( \W_\Phi ) = \sum_{k=0}^n \sum_{(([k])|([n-k])) \wedge \Phi = \Psi|\Gamma}
\W_\Psi \otimes \W_\Gamma.
\end{equation}

\begin{remark} It is interesting to remark that one can define on $\NCQ$ a second internal comultiplication corresponding to the substitution $\M_{\Phi}[X]\mapsto \M_{\Phi}[X Y]$. Here we use the lexicographic order to order the alphabet $X Y$. This gives 
\begin{equation}
\Delta^{\odot}( \W_\Phi ) = \sum_{\Psi \wedge \Gamma = \Psi}
\W_\Psi \otimes \W_\Gamma.
\end{equation}
The dual to this operation gives an algebra on $\NCQ^\ast_n$ for all $n$.
This is precisely the dual of the Solomon-Tits algebra initially defined by 
Tits in \cite{T} (see also \cite{BHRZ}).
\end{remark}

Note that because there is an embedding of $\theta : \NCSI \rightarrow \NCQ$
then the dual of this map is a projection from $\NCQ^\ast$ to $\NCSI^\ast$
which may be given explicitly as the surjection 
$\theta^\ast: \NCQ^\ast \rightarrow \NCSI^\ast$ by $\theta^\ast( \W_\Phi) = \w_{A(\Phi)}$.
This of course follows from the defining relation $[ \theta( \m_A ), \W_\Phi ]_{\NCQ} =
[ \m_A, \theta^\ast( \W_\Phi ) ]_{\NCSI}$.

\end{subsection}
\end{section}

\begin{section}{ $\NCQ$ is free } \label{sec:ncqsymfree}
There is another  order on set compositions defined with the 
covering relations $\Phi \lessdot (\Phi_1, \ldots, \Phi_i \cup \Phi_{i+1}, \Phi_{i+2}, \ldots,
\Phi_{\ell(\Phi)})$ for each $1 \leq i < \ell(\Phi)$ such that every integer in $\Phi_i$
is less than every integer in $\Phi_{i+1}$.  This new order is analogous to the
second order that we defined on set partitions and so we use the
notation $\mleq$ to denote the closure of these covering relations.

Under this order the set composition
$(n\vbar n-1\vbar\cdots\vbar1)$ is not comparable to any other set composition.
In general,
the connected component to any particular set composition will be a boolean lattice.
This is easy to see since any set composition is greater than or equal to
a minimal element $\Pi$ which has all of the number in $[n]$ in separate parts (essentially,
a permutation).
It is easy to see that the elements which are above each of the permutations form
a boolean lattice on this order which is isomorphic to the 
boolean lattice of the set of $\{i \,|\, \Pi_i < \Pi_{i+1}\}$.  In the diagram
below we see how the set compositions of size $3$ are simply the union of boolean
lattices.

\begin{figure}[htbp] 
   \centering
   \includegraphics[width=11cm]{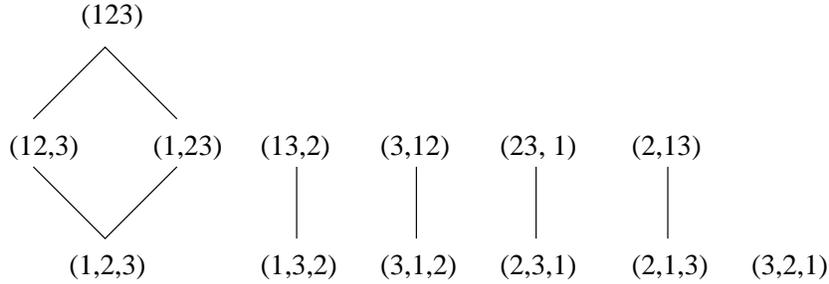} 
   \caption{Hasse diagram of $(\{ \Phi \models [3] \}, \geq_\ast)$}
   \label{fig:example}
\end{figure}
\setlength{\unitlength}{0.00083333in}

Now we can define a basis $\{ \Q_\Phi \}_\Phi$ using this order by
\begin{equation}
\Q_\Phi = \sum_{\Phi' \geq_\ast \Phi} \M_{\Phi'}.
\end{equation}

This basis is not multiplicative as was the $\q_A$ basis of $\NCSI$, but it does
have the following elegant rule for the product of two basis elements. For two set compositions $\Phi\models[n]$ and $\Psi\models[m]$, we define $\Phi\widetilde{\shuf}\Psi$ to be the (non-commutative) shuffle $\Phi\shuf\Psi\!\uparrow_{n}$ where $\Phi$ is viewed as a word in the subset of $[n]$ and $\Psi\!\uparrow_{n}$ is viewed as the word in the subset of $\set{n+1,n+2,\ldots, n+m}$ obtained from $\Psi$ by adding $n$ to every entries.

\begin{theorem}  For $\Phi \models [n]$ and $\Psi \models [m]$,
$$\Q_\Phi \Q_\Psi = \sum_{\Gamma \in \Phi \widetilde{\shuf} \Psi} \Q_\Gamma$$
\end{theorem}

\begin{proof}
From the definitions we have
 $$\Q_\Phi \Q_\Psi =\sum_{\Phi'\geq_\ast \Phi \atop \Psi'\geq_\ast \Psi} \M_{\Phi'}\M_{\Psi'}
          =\sum_{\Phi'\geq_\ast \Phi \atop \Psi'\geq_\ast \Psi} \sum_{(([n])|([m]))\wedge \Gamma'=\Phi'|\Psi'}
                      \M_{\Gamma'}\,.
 $$
On the other hand we have
 $$\sum_{\Gamma \in \Phi \widetilde{\shuf} \Psi} \Q_\Gamma =
      \sum_{\Gamma \in \Phi \widetilde{\shuf} \Psi} \sum_{\Gamma'\geq_\ast \Gamma} \M_{\Gamma'}.
      $$
We first remark that both equations are multiplicity free. In the right hand side of the second equation, we have that $\Gamma'$ satisfies $\Gamma'\geq_\ast\Gamma$ for $\Gamma\in \Phi \widetilde{\shuf} \Psi$. This happen if and only if $(\Phi|\Psi)\wedge\Gamma'=\Phi|\Psi$ which is equivalent to $(([n])|([m]))\wedge \Gamma'=\Phi'|\Psi' \geq_\ast\Phi|\Psi$. This gives us the terms in the right hand side of the first equation and conclude the desired equality.
\end{proof}

For example,
\begin{equation*}
\Q_{(13\vbar2)} \Q_{(1\vbar2)} = 
\Q_{(13\vbar2\vbar4\vbar5)} + \Q_{(13\vbar4\vbar2\vbar5)} +  \Q_{(13\vbar4\vbar5\vbar2)} +
\Q_{(4\vbar13\vbar2\vbar5)} + \Q_{(4\vbar13\vbar5\vbar2)} +  \Q_{(4\vbar5\vbar13\vbar2)}.
\end{equation*}
Notice that this is also a non-commutative product even though we are using the shuffle
operation.

This result then shows the consequence which is one of the main goals
of this paper, namely,
\begin{theorem}
The algebra $\NCQ$ is freely generated by the elements
$$\{ \Q_\Phi \,|\, \Phi \models [n] , n\geq 0, \Phi \hbox{ is atomic}\}.$$
\end{theorem}

\begin{proof}
For $\Phi=(\Phi_1,\Phi_2,\ldots,\Phi_\ell)$ we define $w(\Phi)=w_1 w_2 \cdots w_n$ 
to be the word such that $w_i =j$ if and only if $i\in\Phi_j$. 
For two set composition $\Phi$ and $\Psi$, we say that 
$\Phi\preceq_{lex}\Psi$ if $w(\Phi)\leq_{lex}w(\Psi)$.

We have that if  $\Gamma\in\Phi^{(1)}\widetilde{\shuf}\Phi^{(2)}\widetilde{\shuf}\cdots\widetilde{\shuf}\Phi^{(k)}$,
then $\Phi^{(1)}|\Phi^{(2)}|\cdots|\Phi^{(k)} \preceq_{lex}\Gamma$. This implies that if $\Phi=\Phi^{(1)}|\Phi^{(2)}|\cdots|\Phi^{(k)}$ is the unique minimal atomic decomposition of $\Phi$, then
  $$\Q_{\Phi^{(1)}}\Q_{\Phi^{(2)}}\cdots \Q_{\Phi^{(k)}} = \Q_\Phi + \sum_{\Phi\prec_{lex}\Gamma} c_\Gamma \Q_\Gamma.$$
The theorem follows by triangularity.
\end{proof}

It is also interesting that because the order that we have chosen is an
Eulerian poset, we know that the dual elements in $\NCQ^\ast$
are easy to calculate.  In fact, we have
$$\Q_\Phi^\ast = \sum_{\Psi \geq_\ast \Phi} (-1)^{\ell(\Phi) - \ell(\Psi)} \W_\Psi.$$

\end{section}

\begin{section}{ $\NCQ$ is cofree }\label{sec:ncqsymcofree}
So that we can again establish a basis which has elegant
multiplicative properties in the algebra of $\NCQ^\ast$, we
need to introduce another order on set compositions.  The
orders which we define for these bases arise because they
interact nicely with respect to the product operations on the
algebras, what is remarkable is that at the same time the
posets have properties which make them interesting to
study in their own right. 
In a related paper with C. Hohlweg \cite{BHZ} we consider some
of the properties of this poset in more detail.

We will
define a  order in which only elements which have have the same
image under $\alpha$ are comparable. 
Consider the order which we define $\Phi \leq_\# \Psi$ if
$\alpha(\Phi) = \alpha(\Psi)$ and if $\Phi^! = ( \Phi^{(1)}, \Phi^{(2)}, \ldots, \Phi^{(k)})$
then $\Phi^{(i)} = st( \Psi_{\ell_{i-1}+1}, \Psi_{\ell_{i-1}+2}, \ldots, \Psi_{\ell_{i}})$
where $\ell_0 = 0$ and $\ell_i = \ell(\Phi^{(1)})+\ell(\Phi^{(2)})+\cdots + \ell(\Phi^{(i)})$.  
In the appendix we include the Hasse diagrams for this poset for the set compositions
of size $3$ and $4$.

For example, $(123,4,56,7)$ is smaller than or equal to all set compositions $\Psi$ such
that $\alpha(\Psi) = (3,1,2,1)$ since $(123,4,56,7)^! = ((123),(1),(12),(1))$ and 
$st(\Psi_1) = (123)$, $st(\Psi_2) = (1)$, $st(\Psi_3) = (12)$ and $st(\Psi_4) = (1)$ 
for all $\Psi \models [7]$ with
$\alpha(\Psi) = (3,1,2,1)$.  $(134,2,56,7)$ is smaller than $(156,2,37,4)$ because
$(134,2,56,7)^! = ((134,2), (12), (1)) = (st(156,2), st(37), st(1))$, but it is
not comparable to say $(124,5,36,7)$ since $(st(124,5), st(36), st(7)) = ((123,4), (12), (1))$.

In the appendix we give drawings of the poset for $n=4$ representing the
elements of this  order.  In particular, we note that for the set compositions
$\Phi  \models [n]$ with $\alpha(\Phi) = (1^n)$ this defines a  order on
permutations (since there is a natural correspondence between permutations and
set compositions with $\alpha(\Phi) = (1^n)$) where the number of elements
of rank $k$ are given by the coefficient of $x^n t^k$ in the generating function
$\frac{g(x)}{t+g(x)-t g(x)}$ where $g(x) = \sum_{k\geq 0} k! x^k$ 
(\cite{OLEIS} sequence A059438). See \cite{BHZ} for more results 
about this order and other related ones.

This order appears naturally in the algebra of $\NCQ^\ast$ because if we
define a basis
$$\V_\Phi = \sum_{\Phi' \geq_\# \Phi} \W_{\Phi'}$$
then the basis has very elegant rule for computing the product.
\begin{theorem} For $\Phi \models [n]$ and $\Psi \models [k]$, then
\begin{equation*}\label{cofree2}
\V_\Phi \V_\Psi = \V_{\Phi|\Psi}
\end{equation*}
\end{theorem}

\begin{proof}
We have
 $$\V_\Phi \V_\Psi =\sum_{\Phi'\geq_\# \Phi \atop \Psi'\geq_\# \Psi} \W_{\Phi'}\W_{\Psi'}
          =\sum_{\Phi'\geq_\# \Phi \atop \Psi'\geq_\# \Psi} \sum_{S\in {[n+m] \choose n}}
                      \W_{\Phi'\!\uparrow_S \cdot \Psi'\!\uparrow_{S^c}}\,.
 $$
On the other hand we have
 $$\V_{\Phi|\Psi} =
      \sum_{\Gamma \geq_\# \Phi|\Psi} \W_{\Gamma}.
      $$
We remark that both equations are multiplicity free. 
For any $S\in {[n+m] \choose n}$, $\Phi'\geq_\# \Phi$ and  $\Psi'\geq_\# \Psi$ we have that
  $$\alpha\big(\Phi'\!\uparrow_S \cdot \Psi'\!\uparrow_{S^c}\big) = \alpha(\Phi'|\Psi')
        =\alpha(\Phi')\cdot\alpha(\Psi')=\alpha(\Phi)\cdot\alpha(\Psi)=\alpha(\Phi|\Psi).
        $$
Also we have that $(\Phi|\Psi)^!=\Phi^!\cdot\Psi^!=(\Phi^{(1)},\ldots,\Phi^{(k)},\Psi^{(1)},\ldots,\Psi^{(r)})$. Since
 $$\Phi^{(i)}=st\big(\Phi'_{\ell_{i-1}+1},\ldots,\Phi'_{\ell_{i}}\big)=st\big(\Phi'_{\ell_{i-1}+1}\!\uparrow_S,\ldots,\Phi'_{\ell_{i}}\!\uparrow_S\big)$$
 and
  $$\Psi^{(j)}=st\big(\Psi'_{t_{j-1}+1},\ldots,\Phi'_{t_j}\big)=st\big(\Phi'_{t_{j-1}+1}\!\uparrow_S,\ldots,\Phi'_{t_{j}}\!\uparrow_S\big)$$
   for $\ell_i=\ell(\Phi^{(1)})+\cdots +\ell(\Phi^{(i)})$ and $t_j=\ell(\Psi^{(1)})+\cdots +\ell(\Psi^{(j)})$, we have that 
$\Phi'\!\uparrow_S \cdot \Psi'\!\uparrow_{S^c}\geq_\#\Phi|\Psi$. Conversely, if $\Gamma\geq_\#\Phi|\Psi$, there is a unique $S$ such that $\Gamma=\Phi'\!\uparrow_S \cdot \Psi'\!\uparrow_{S^c}$. This shows the desired equality.
\end{proof}

For example if we compute $\V_{(12)} = \W_{(12)}$ and $\V_{(1\vbar2)} = \W_{(1\vbar2)}
+ \W_{(2\vbar1)}$ then 
\begin{align*}
\V_{(12)} \V_{(2\vbar1)} &= \W_{(12)}\W_{(2\vbar1)}\\
&= 
\W_{(12\vbar4\vbar3)}
+ \W_{(13\vbar4\vbar2)}+ \W_{(14\vbar3\vbar2)}+ \W_{(23\vbar4\vbar1)}
+ \W_{(24\vbar3\vbar1)}+ \W_{(34\vbar2\vbar1)}\\
&= \V_{(12,4,3)}.~
\end{align*}

In light of Theorem~\ref{cofree2} it is clear that the set $\set{\V_\Phi\,|\,\Phi\models[n], n\geq 0, \Phi \hbox{ is atomic}  }$ freely generate $\NCQ^\ast$. This implies our last  theorem.
\begin{theorem} The algebra $\NCQ$ is cofree.
\end{theorem}

\label{numsections}
\end{section}

\begin{section}{Appendix A: Poset of $(\{ \Phi \models [4]\}, \leq_{\#})$ and 
$(\{ \Phi \models [3]\}, \leq_{\#})$}

For a set composition $\Phi$, let $\alpha(\Phi)$ be the composition
$(|\Phi_1|, |\Phi_2|, \ldots, |\Phi_{\ell(\Phi)}|)$.  The poset of
$(\{ \Phi \models [n] | \alpha(\Phi) = (\alpha_1, \alpha_2, \ldots, \alpha_k)\}, \leq_{\#})$ is
isomorphic to the poset 
$(\{ \Phi \models [n] | \alpha(\Phi) = (\alpha_k, \alpha_{k-1}, \ldots, \alpha_1)\}, \leq_{\#})$
by reversing the entries in $\Phi$ and complementing the entries.  This appendix
includes the Hasse diagrams for the set compositions of size $3$ and $4$.

\begin{figure}[htbp] 
   \centering 
   \includegraphics[width=12.4cm]{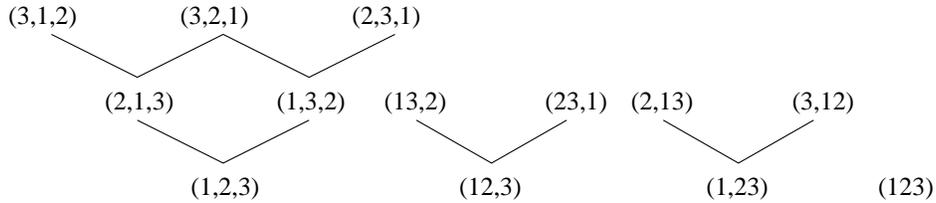} 
   \caption{Poset of $\{ \Phi \models [3]  \}$}
\end{figure}

\begin{figure}[htbp] 
   \centering 
   \includegraphics[width=4.6cm]{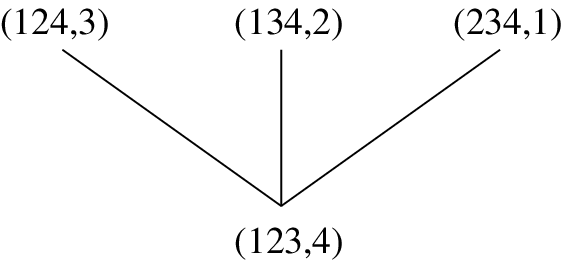} 
   \caption{Poset of $\{ \Phi \models [4] : \alpha(\Phi) = (3,1) \}$}
   \label{fig:example}
\end{figure}

\begin{figure}[htbp] 
   \centering 
   \includegraphics[width=8.1cm]{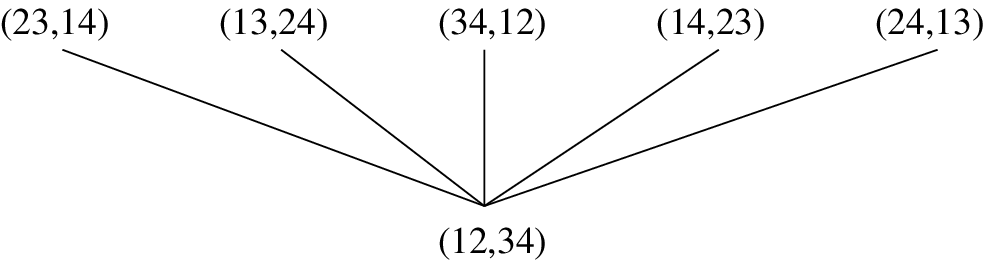} 
   \caption{Poset of $\{ \Phi \models [4] : \alpha(\Phi) = (2,2) \}$}
   \label{fig:example}
\end{figure}

\begin{figure}[htbp] 
   \centering 
   \includegraphics[width=12.0cm]{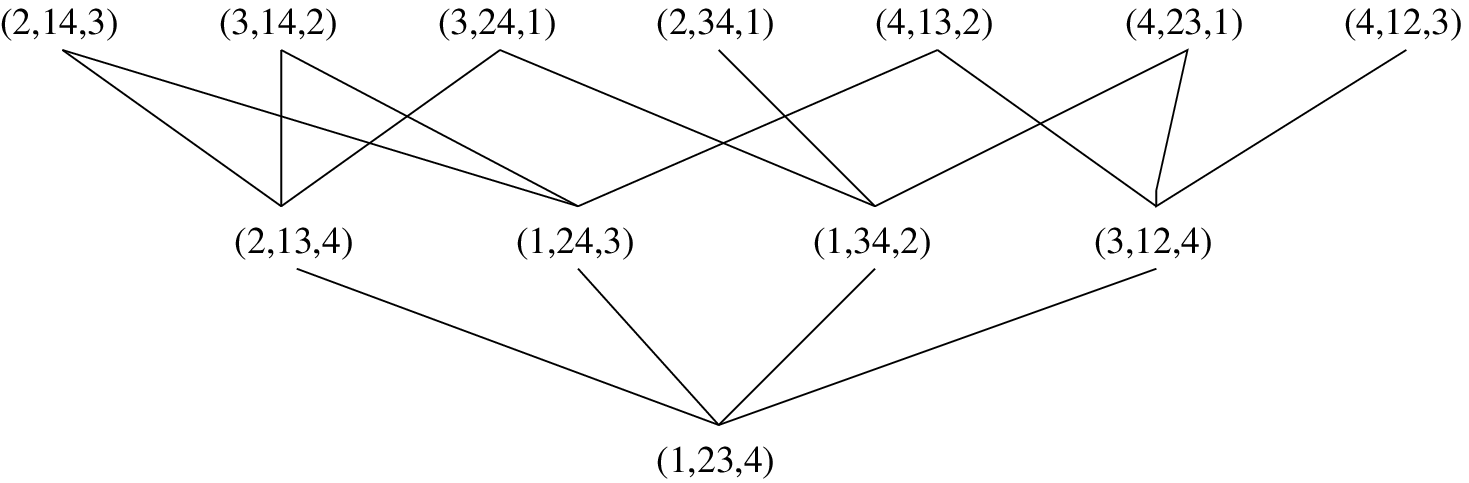} 
   \caption{Poset of $\{ \Phi \models [4] : \alpha(\Phi) = (1,2,1) \}$}
   \label{fig:example}
\end{figure}

\begin{figure}[htbp] 
   \centering 
   \includegraphics[width=11.7cm]{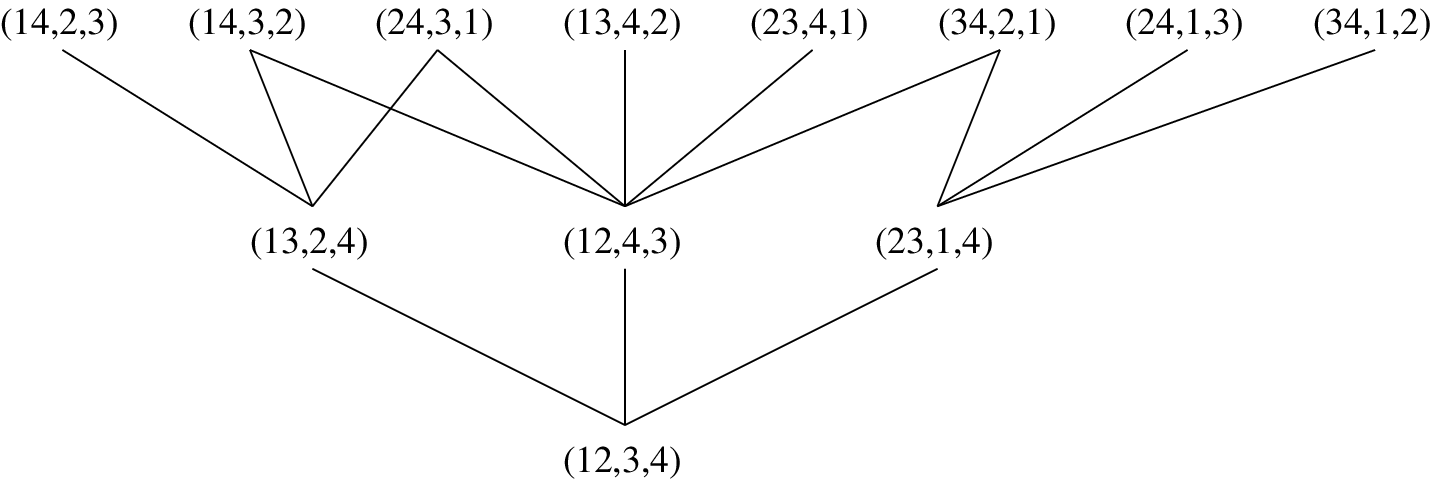} 
   \caption{Poset of $\{ \Phi \models [4] : \alpha(\Phi) = (2,1,1) \}$}
   \label{fig:example}
\end{figure}

\begin{figure}[htbp] 
   \centering
   \includegraphics[width=6.5in]{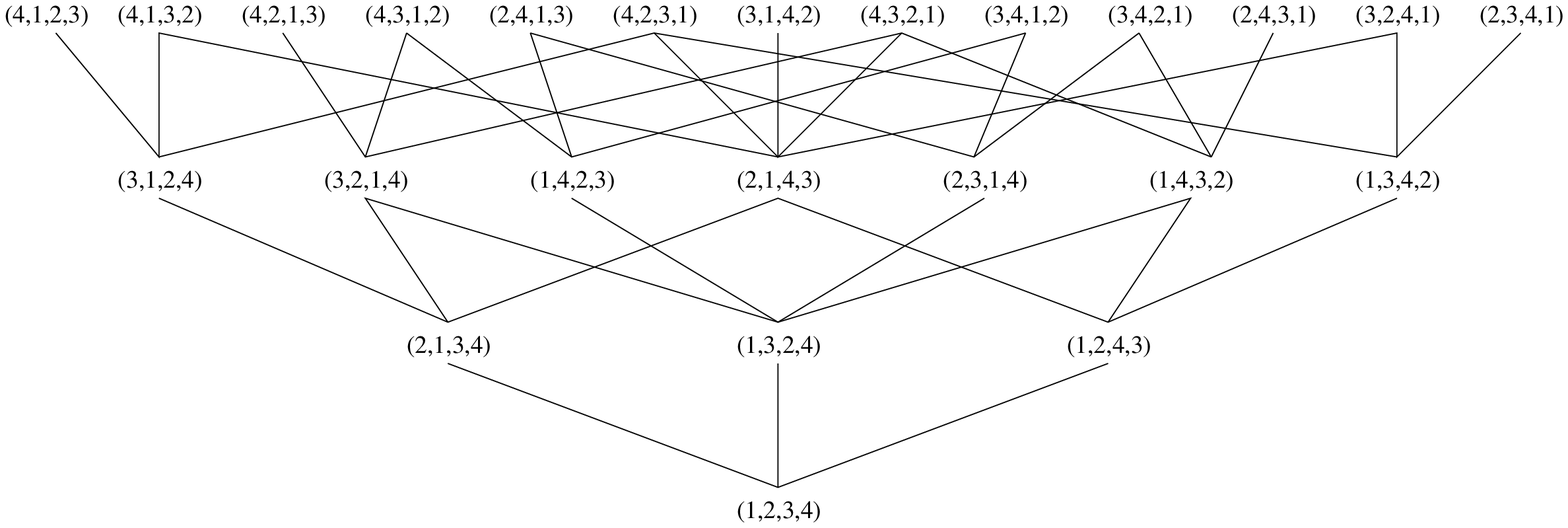} 
   \caption{Poset of $\{ \Phi \models [4] : \alpha(\Phi) = (1,1,1,1) \}$}
   \label{fig:example}
\end{figure}

\end{section}

\end{document}